\documentclass[11pt,
a4]{article}
\usepackage{fancybox,%
            epsf,%
            theorem,%
            eepic,%
            psfig,%
            epsfig,%
            amssymb,%
            amsmath,%
            multicol,%
            pifont,%
}

\usepackage%
[%
nodraft
]%
{symbols}
 \let\ac=\coauthorone                 
 \def\coauthoronesname{AC}
 
\begin{document}

\title{\Large  Uniform Semiglobal Practical Asymptotic Stability for Non-autonomous Cascaded Systems and Applications}
\author{\ \\[-3mm]  Antoine Chaillet \qquad Antonio Lor\'{\i}a\\[1mm] 
\address{
CNRS--LSS, Sup\'elec, 3 rue Joliot Curie,  91192 Gif s/Yvette, France, {\tt chaillet@lss.supelec.fr}, {\tt loria@lss.supelec.fr} \\
}}

 \date{ \small \today \ \  \em 
 }

 \sloppy
 \maketitle
 \parskip = 4pt

\begin{abstract}
It is due to the modularity of the analysis that results for cascaded systems have proved their utility in numerous control applications as well as in the development of general control techniques based on ``adding integrators''. Nevertheless, the standing assumptions in most of the present literature on cascaded systems is that, when \emph{decoupled}, the subsystems constituting the cascade are uniformly globally asymptotically stable (UGAS). Hence existing results fail in the more general case when the subsystems are uniformly semiglobally practically asymptotically stable (USPAS). This situation is often encountered in control practice, e.g., in control of physical systems with external perturbations, measurement noise, unmodelled dynamics, etc.  This paper generalizes previous results for cascades by establishing that, under a uniform boundedness condition, the cascade of two USPAS systems remains USPAS. An analogous result can be derived for USAS systems in cascade. Furthermore, we show the utility of our results in the PID control of mechanical systems affected by unknown non-dissipative forces and considering the dynamics of the DC motors.
\end{abstract}

\section{Introduction}
Cascaded dynamical systems appear in many control applications whether naturally or intentionally due to control design. Cascades-based control consists in designing the control law so that the closed loop system has a cascaded structure. Such strategy has often the advantage of reducing the complexity of the controller and the difficulty of the stability analysis (see e.g. \cite{KKK,ERJENSTHESIS,CASCFAP}) as opposed to more general Lyapunov-based control methods. From a theoretical viewpoint the problem of stability analysis of cascaded systems has attracted the interest of the community starting with the seminal paper \cite{SUSKOK}. See also \cite{SEPBOOK} and references therein. In general terms, a fundamental result that one may retain from the literature is that cascades of uniformly globally asymptotically stable systems (UGAS) remain UGAS if and only if the solutions are uniformly globally bounded (UGB). See \cite{SEISUA,SON89} for the proof of this statement in the case of autonomous systems and \cite{CASCAUT} for the case of time-varying systems. In a similar spirit, see \cite{SONTAC03} for a local result: a proof that a LAS system perturbed by a converging input which is such that its solutions remain in the domain of attraction of the nominal system, remains LAS. Other works on stability of cascades deal, directly or indirectly, with the fundamental problem of establishing conditions for (uniform global) boundedness of the trajectories. A sufficient condition is the well-understood notion of input-to-state stability (ISS). 

A considerable drawback of most results on stability of cascades available in the literature is that they rely on the assumption of {\em global} stability properties for the separate subsystems. Nevertheless, in practice, it is often the case that only local (with a specified estimate of the domain of attraction) or {\em semiglobal} properties can be concluded. Semiglobal asymptotic stability pertains to the case when one can prove that, by tuning certain parameter of the control system, the estimate of the domain of attraction can be arbitrarily enlarged. Such parameter is often, but not always, a control gain. We are not aware of results on semiglobal asymptotic stability nor ISS for cascades of ``semiglobal ISS'' systems. Another case which is not covered by most results in the literature of cascaded systems is that of stability with respect to balls. A notable exception, in a way, is \cite{JIATEEPRA} which introduced the notion of input to state {\em practical} stability. However, this notion is different from the one used here in the sense that the neighborhood of the origin which is stable is not required to be arbitrarily reducible by a convenient choice of a parameter. More generally, we are not aware of results on semiglobal practical asymptotic stability of cascades for {\em continuous} time-varying systems. For (parameterized) \emph{discrete-time} systems the only results that we know of are those presented in \cite{DTCASCTAC}. Needless to say that the nature of discrete-time non-autonomous systems parameterized in the sampling time is fundamentally different.

In this paper we address the stability analysis problem for cascades of time-varying systems that are uniformly semiglobally practically asymptotically stable (USPAS). We establish that, under a uniform boundedness condition on its solutions and provided a Lyapunov function for the ``perturbed" subsystem, the cascade of two USPAS systems remains USPAS. In the same way, the cascade of two USAS systems is USAS. Our main result extends in this direction,  \cite[Lemma 2]{CASCAUT} and the main results of \cite{SEISUA,SON89} which are, in a way, at the basis of many theorems on UGAS of cascades. 

The rest of the paper is organized as follows. In next section we present some definitions of stability and an auxiliary proposition on semiglobal practical asymptotic stability. Our main result is presented in Section \ref{sec:main}. In Section \ref{sec:eg} we illustrate the utility of our findings with an example dealing with the PID control of a mechanical system taking into account the dynamics of the DC motors. The proofs of all the results are given in Section \ref{sec:proofs} and we conclude with some remarks in Section \ref{sec:concl}.

\section{Definitions and preliminary results}

\noindent {\bf Notation.}
A continuous function $\alpha : \mRp \to \mRp$ is of \emph{class ${\cal K}$} ($\alpha \in {\cal K}$), if it is strictly increasing and $\alpha(0) = 0$; $\alpha \in {\cal K}_\infty$ if, in addition, $\alpha(s) \to \infty$ as $s\to \infty$. A continuous function $\sigma: \mRp\to\mRp$ is of class $\cal L$ ($\sigma\in\cal L$) if it is non-increasing and tends to zero as its argument tends to infinity. A function $\beta: \mRp\times \mRp \to \mRp$ is said to be a \emph{class $\cal{KL}$ function} if, $\beta(\cdot,t)\in\cal K$ for any $t\geq 0$, and $\beta(s,\cdot)\in\cal L$ for any $s\geq 0$. We denote by $x(\cdot,t_0,x_0)$ the solutions of the differential equation $\dot x = f(t,x)$ with initial conditions $(t_0,x_0)$. We use $\norm{\cdot}$ for the Euclidean norm of vectors and the induced $L_2$ norm of matrices. We denote by $\mathcal B_\delta$ the \emph{closed} ball in $\mathbb R^n$ of radius $\delta$, {\em i.e.} $\mathcal B_\delta:=\{x\in\mR^n~|~\norm{x}\leq \delta\}$. We use the notation $\mathcal H(\delta,\Delta):=\left\{x\in\mR^n ~|~ \delta\leq\norm x\leq\Delta\right\}$. We define $\norm{x}_\delta:=\textrm{inf}_{z\in\mathcal B_\delta}\norm{x-z}$. When the context is sufficiently explicit, we may omit to write the arguments of a function. 

\subsection{Asymptotic stability of balls}
\vspace{-0.5cm}For nonlinear time-varying systems of the form
\begin{equation}\label{sys1}
\dot x=f(t,x)\,,
\end{equation}
where $x\in\mR^n$, $t\in\mRp$ and $f:\mRp\times\mR^n \to \mR^n$ is piecewise continuous in $t$ and locally Lipschitz in $x$, we introduce the following.
\begin{defin}[US of a ball]\label{defULS}
Let $\delta$ and $\Delta$ be nonnegative numbers such that $\Delta>\delta$. The ball $\mathcal B_\delta$ is said to be \emph{Uniformly Stable on $\mathcal B_\Delta$} for the system \rref{sys1} if there exists a class $\mathcal{K}_\infty$ function $\eta$ such that the solutions of \rref{sys1} from any initial state $x_0 \in \mathcal B_\Delta$ and initial time $t_0\in\mRp$ satisfy
$$
\norm{x(t,t_0,x_0)}_{\delta} \leq \eta(\norm{x_0})\,, \quad \forall t\geq t_0\,.
$$

\end{defin}
\begin{defin}[UA of a ball]\label{defULA}
Let $\delta$ and $\Delta$ be nonnegative numbers such that $\Delta>\delta$. The ball $\mathcal B_\delta$ is said to be \emph{Uniformly Attractive on $\mathcal B_\Delta$} for the system \rref{sys1} if there exists a class $\mathcal{L}$ function $\sigma$ such that the solutions of \rref{sys1} from any initial state $x_0 \in \mathcal B_\Delta$ and initial time $t_0\in\mRp$ satisfy
\begin{displaymath}
\norm{x(t,t_0,x_0)}_{\delta} \leq \sigma(t-t_0)\,, \quad \forall t\geq t_0\,.
\end{displaymath}
\end{defin}
\begin{defin}[UAS of a ball]\label{defULAS}
Let $\delta$ and $\Delta$ be nonnegative numbers such that $\Delta>\delta$. The ball $\mathcal B_\delta$ is said to be \emph{Uniformly Asymptotically Stable on $\mathcal B_\Delta$} for the system \rref{sys1} if it is both US and UA on $\mathcal B_\Delta$. 
\end{defin}

\begin{rmk}\label{rmkstrongULAS}
The property of UAS of a ball defined above is less restrictive than the time-varying adaptation of ``asymptotic stability with respect to a set'' given in \cite{LINSONWAN} for the case when the set is a ball. Indeed, in the latter reference, it is imposed that the ball $\mathcal B_\delta$ be positively invariant. Notice also that, modulo that ${\cal B}_\Delta$ and ${\cal B}_\delta$ are closed, Definition \ref{defULAS} is equivalent to $\Delta\to\delta$ stability as defined in \cite{TEPEAE98}; more precisely, UAS of $\mathcal B_\delta$ on $\mathcal B_\Delta$ implies $\Delta\to\delta$ stability and $\Delta\to\delta$ stability implies UAS of $\mathcal B_{\delta}$ on $\mathcal B_{\Delta'}$, for all $\Delta'\in(\delta,\Delta)$.
\end{rmk}
We also emphasize that the UAS of a ball can be characterized by a $\cKL$ bound as in the case of UAS of the origin.
\begin{propos}[$\mathcal{KL}$ estimate for UAS of a ball]\label{propKL}
The ball $\mathcal B_\delta$ is UAS on $\mathcal B_\Delta$ for the system \rref{sys1} if and only if there exists a class $\mathcal{KL}$ function $\beta$ such that the solutions of \rref{sys1} from any initial state $x_0 \in \mathcal B_\Delta$ and initial time $t_0\in\mRp$ satisfy
\begin{displaymath}
\norm{x(t,t_0,x_0)}_{\delta} \leq \beta(\norm{x_0},t-t_0)\,, \quad \forall t\geq t_0\,.
\end{displaymath}
\end{propos}

\noindent Our main result for cascades is formulated based on the following property of boundedness of solutions.

\begin{defin}[UB] \label{defULB} 
The solutions of \rref{sys1} are said to be \emph{Uniformly Bounded on the compact set $\mathcal A\subset\mR^n$} if there exist a class $\mathcal K$ function $\gamma$ and a nonnegative constant $\mu$ such that, for any initial time $t_0\in\mRp$ and any initial state $x_0\in\mathcal A$, it holds that
\begin{displaymath}
\norm{x(t,t_0,x_0)}\leq \gamma(\norm{x_0})+\mu\,,\quad \forall t\geq t_0\,.
\end{displaymath}
\end{defin}

\subsection{Semiglobal practical asymptotic stability} 
\vspace{-0.5cm}Our main result addresses the problem of uniform semiglobal practical asymptotic stability (USPAS) for parameterized nonlinear time-varying systems of the form
\begin{equation}\label{sys2}
\dot x=f(t,x,\theta)\,,
\end{equation}
where $x\in\mR^n$, $t\in\mRp$, $\theta\in\mR^m$ is a constant parameter and $f:\mRp\times\mR^n\times\mR^m \to \mR^n$ is locally Lipschitz in $x$ and piecewise continuous in $t$.

\begin{defin}[USPAS] \label{defUSPAS}
Let $\Theta\subset\mR^m$ be a set of parameters. The system \rref{sys2} is said to be \emph{Uniformly Semiglobally Practically Asymptotically Stable on  $\Theta$} if, given any $\Delta>\delta >0$, there exists $\theta^\star\in\Theta$ such that $\mathcal B_\delta$ is UAS on $\mathcal B_\Delta$ for the system $\dot x=f(t,x,\theta^\star)$. 
\end{defin}

In view of Proposition \ref{propKL}, the above definition is equivalent to the following statement: for any given pair $\Delta>\delta>0$, there exists a parameter $\theta^\star\in\Theta$ and a $\cKL$ function $\beta$ such that, for all $x_0\in\mathcal B_\Delta$ and all $t_0\in\mRp$, $\norm{x(t,t_0,x_0,\theta^\star)}_\delta\leq \beta(\norm{x_0},t-t_0)$ for all $t\geq t_0$. We stress that, although this does not explicitly appear in the notation, the function $\beta$ is not required to be independent of $\delta$ and $\Delta$, as opposed to some other definitions of semiglobal and/or practical stability existing in the literature, such as \cite{TEPEAE98,DTCASCTAC}. This non-uniformity in $\delta$ and $\Delta$ makes our definition of USPAS a more general property since, as it will appear more clear in the sequel, it allows to make use of Lyapunov functions with bounds that may depend on the tuning parameter $\theta$.

Systems of the form \rref{sys2} result, for instance, from control systems in {\em closed loop}; in this case, we can think of $\theta$ as control gains or other design parameters. Then we say that \rref{sys2} is USPAS if the estimate of the domain of attraction $\mathcal B_\Delta$ and the ball $\mathcal B_\delta$ which is UAS can be arbitrarily enlarged and diminished, respectively, by a convenient choice of the design parameters. Such situation is fairly common in control practice; for instance, in output feedback tracking control of mechanical systems ({\it cf.} \cite{TEEPRA1}). For discrete-time systems one also finds it useful to define USPAS with respect to a design parameter: in this case, the sampling time. See \cite{NESTEEPK,NESTEESON} for definitions and a solid framework on USPAS for discrete-time systems. 

Notice also that when, by an abuse of notation, $\delta=0$ we recover from Definition \ref{defUSPAS} the notion of uniform semiglobal asymptotic stability (USAS). If, in addition, $\Delta=\infty$, we recover the definition of uniform global asymptotic stability (UGAS).

The following result gives a sufficient condition, in terms of a Lyapunov function, for the dynamical parameterized system \rref{sys2} to be uniformly semiglobally practically asymptotically stable on a given set of parameters. See Section \ref{proofprop1} for the proof.

\begin{propos}[Lyapunov sufficient condition for USPAS]\label{prop1} 
Suppose that, given any $\Delta>\delta>0$, there exist a parameter $\theta\in \Theta$, a continuously differentiable Lyapunov function\footnote{It should be clear that $V$ may depend on $\delta$ and $\Delta$ as well. We shorten the notation $V_{\delta,\Delta}$ to just $V$ for clarity.} $V:\mRp\times\mR^n\to\mRp$, class $\mathcal K_\infty$ functions $\underline\alpha_{\delta,\Delta}$, $\overline\alpha_{\delta,\Delta}$, $\alpha_{\delta,\Delta}$, and a continuous positive nondecreasing function $c_{\delta,\Delta}$ such that, for all $x\in\mathcal H(\delta,\Delta)$ and all $t\in\mR_{\geq 0}$,
\begin{equation}\label{pilot1}
\underline\alpha_{\delta,\Delta}(\norm {x}_{\delta}) \leq V(t,x) \leq \overline\alpha_{\delta,\Delta}(\norm {x})
\end{equation}
\begin{equation}\label{pilot2}
\frac{\partial V}{\partial t}(t,x)+\frac{\partial V}{\partial x}(t,x)f(t,x,\theta) \leq -\alpha_{\delta,\Delta}(\norm x)
\end{equation}
Assume further that, for any $\Delta>0$,
\begin{equation}\label{condadd}
\lim_{\delta\to 0} \underline\alpha_{\delta,\Delta}^{-1}\circ \overline\alpha_{\delta,\Delta}(\delta)=0
\end{equation}
and that, for any $\delta>0$,
\begin{equation}\label{condadd2}
\lim_{\Delta\to \infty} \overline\alpha^{-1}_{\delta,\Delta}\circ \underline\alpha_{\delta,\Delta}(\Delta)=\infty\,.
\end{equation}
Then, the system \rref{sys2} is USPAS on the parameter set $\Theta$.
\end{propos}

It is worth mentioning that the condition $\underline\alpha_{\delta,\Delta}(\norm {x}) \leq V(t,x) \leq \overline\alpha_{\delta,\Delta}(\norm {x})$, which \emph{implies} \rref{pilot1}, often holds in the analysis of control systems. In particular, it holds for systems with additive bounded disturbances when USPAS may be inferred using a Lyapunov function for UGAS of the corresponding unperturbed system. See Section \ref{sec:eg} for an example. Condition \rref{pilot2} also appears naturally in the context of stability of perturbed systems. The last two conditions, \rref{condadd} and \rref{condadd2}, need to be imposed due to the fact that the Lyapunov function $V$ is {\em not} required to be the same for all $\delta$ and all $\Delta$. As we show in the proof, conditions \rref{condadd} and \rref{condadd2} ``compensate'' this  lack of uniformity to ensure that the estimate of the domain of attraction $\mathcal B_\Delta$ and the set $\mathcal B_\delta$ which is UAS can be arbitrarily enlarged and diminished respectively.

\begin{rmk}\label{rmkprop1}
By noticing that the UAS of $\mathcal B_\delta$ on $\mathcal B_\Delta$ implies the UAS of $\mathcal B_{\delta'}$ on $\mathcal B_{\Delta'}$ for any $\delta'$ and $\Delta'$ satisfying $\delta\leq \delta'<\Delta'\leq \Delta$, the conclusion of Proposition \ref{prop1} remains valid if \rref{pilot1} and \rref{pilot2} hold for all $\delta$ \emph{small enough} and all $\Delta$ \emph{large enough}. This relaxed assumption may be useful in practice.
\end{rmk}


\section{Stability of cascades}
\label{sec:main}

We now consider cascaded systems of the form
\begin{subequations}  \label{syscorol}
  \begin{eqnarray}
  \label{syscorol:a}
\dot x_1 & = & f_1(t,x_1,\theta_1) + g(t,x,\theta)x_2 \\ 
  \label{syscorol:b}
\dot x_2 & = & f_2(t,x_2,\theta_2) \,
\end{eqnarray}
\end{subequations}
where $x:=(x_1,x_2)\in\mR^{n_1}\times\mR^{n_2}$, $\theta:=(\theta_1,\theta_2)\in\mR^{n_1}\times\mR^{m_2}$, $t\in\mRp$, $f_1$, $f_2$ and $g$ are locally Lipschitz in the state and piecewise continuous in the time. In order to simplify the statement of our main result, we first introduce the following notation.
\begin{defin}[$\mathcal D$-set]\label{Dset}
For any $\Delta>\delta\geq 0$, the \emph{$\mathcal D$-set} of \rref{sys2} is defined as
\begin{displaymath}
\mathcal D_f(\delta,\Delta):=\left\{\theta \in \mathbb R^m~:~ \mathcal B_\delta \textrm{ is UAS on } \mathcal B_\Delta \textrm{ for } \rref{sys2}  \right\}\,.
\end{displaymath}
\end{defin}

\begin{thm}\label{theo01} 
Under Assumptions \ref{A2}--\ref{A4} below, the cascaded system \rref{syscorol} is USPAS on $\Theta_1\times\Theta_2$.
\end{thm}
\begin{ass}\label{A2}
The system \rref{syscorol:b} is USPAS on $\Theta_2$.
\end{ass}
\begin{ass}[Lyapunov USPAS of the $x_1$-subsystem]\label{A1}
Given any $\Delta_1>\delta_1>0$, there exist a parameter $\theta_1^\star(\delta_1,\Delta_1)\in \Theta_1$, a continuously differentiable Lyapunov function $V_1$, class $\mathcal K_\infty$ functions $\underline\alpha_{\delta_1,\Delta_1}$, $\overline\alpha_{\delta_1,\Delta_1}$, $\alpha_{\delta_1,\Delta_1}$ and a continuous positive nondecreasing function $c_{\delta_1,\Delta_1}$ such that, for all $x_1\in\mathcal H(\delta_1,\Delta_1)$ and all $t\in\mR_{\geq 0}$, bounds \rref{pilot1}, \rref{pilot2} and 
\begin{equation}\label{pilot3}
\norm{\frac{\partial V_1}{\partial x_1}(t,x_1)}\leq c_{\delta_1,\Delta_1}(\norm{x_1})
\end{equation}
hold and, for any positive $\Delta_1$,
\begin{equation}\label{condadd3}
\lim_{\delta_1\to 0} \underline\alpha_{\delta_1,\Delta_1}^{-1}\circ \overline\alpha_{\delta_1,\Delta_1}(\delta_1)=0\,.
\end{equation}
\end{ass}
\begin{ass}[Boundedness of the interconnection term]\label{A23}
The function $g$ is uniformly bounded both in time and in parameter, {\em i.e.} there exists a nondecreasing function $G$ such that, for all  $x\in \mathbb R^{n_1}\times\mathbb R^{n_2}$, all $\theta\in\Theta_1\times\Theta_2$ and all $t\in\mR_{\geq 0}$,
$$
\norm{g(t,x,\theta)}\leq G(\norm{x})\,.
$$
\end{ass}
\begin{ass}[Boundedness of solutions]\label{A4}
There exists a positive constant $\Delta_0$ such that, for any given positive numbers $\delta_1$, $\Delta_1$, $\delta_2$, $\Delta_2$, satisfying $\Delta_1>\max\{\delta_1;\Delta_0\}$ and $\Delta_2>\delta_2$, and for the parameter $\theta_1^\star(\delta_1,\Delta_1)$ as defined in Assumption \ref{A1}, there exists a parameter $\theta_2^\star\in \mathcal D_{f_2}(\delta_2,\Delta_2)\cap\Theta_2$ (cf. Definition \ref{Dset}) and a continuous function\footnote{It should be clear that $\gamma$ may depend on $\delta_1$ and $\delta_2$ as well, but we do not write this dependency explicitly in order to lighten the notation.} $\gamma_{\Delta_1,\Delta_2}:\mR_{>0}\times\mR_{>0}\to\mRp$ such that 
\begin{equation}\label{condgamma}
\lim_{\Delta_1,\Delta_2\to\infty} \gamma_{\Delta_1,\Delta_2}(\Delta_1,\Delta_2) = +\infty\,,
\end{equation}
and the trajectories of \rref{syscorol} with $\theta=\theta^\star$ satisfy
$$
\norm{x_0}\leq \gamma_{\Delta_1,\Delta_2}(\Delta_1,\Delta_2)\qquad\Rightarrow\qquad \norm{x(t,t_0,x_0,\theta^\star)}\leq \Delta_1\,,\qquad \forall t\geq t_0\,.
$$
\end{ass}

The proof of Theorem \ref{theo01} consists in constructing the balls ${\mathcal B_\delta}$ and ${\mathcal B_\Delta}$ and a $\mathcal{KL}$ estimate for the solutions of the cascaded system, based on the respective balls for the $x_1$ ({\em i.e.} \rref{syscorol:a} with $x_2\equiv 0$) and the $x_2$ subsystems. For clarity of exposition we present this in Section \ref{prooftheo01}.

In view of Proposition \ref{prop1}, Assumption \ref{A1} corresponds to the Lyapunov sufficient condition for USPAS of the zero-input $x_1$-subsystem, with the additional condition of a bound on the gradient of $V_1$; we stress that the requirement corresponding to \rref{condadd2} is no longer needed under Assumption \ref{A4}. We state the main result under the more restrictive assumption than simply ``USPAS'' since our proof relies on the explicit knowledge of the Lyapunov function $V_1$. Besides \cite{CHLOautomatica05}, we are not aware of a converse theorem for USPAS, as defined here, that gives all the required properties. In this respect, a converse theorem for USPAS of time-varying systems follows from \cite[Corollary 1]{TEEPRAconverse} for the case of locally Lipschitz $f_1$, but it only establishes the two first inequalities of Assumption \ref{A1} and does not establish that the functions used to bound $V_1$ can be chosen such that \rref{condadd3} holds. Yet, both \rref{pilot3} and \rref{condadd3} are little restrictive, and are satisfied in many concrete applications, as for instance the case study of Section \ref{sec:eg}.

\begin{rmk}\label{rmktheo01}
In view of Remark \ref{rmkprop1}, it is in fact sufficient that the requirements of Assumption \ref{A1} hold for all \emph{small} $\delta$ and all \emph{large} $\Delta$. Also, it is worth pointing out that Assumption \ref{A4} may be relaxed to uniform boundedness on $\mathcal B_{\Delta_1}\times\mathcal B_{\Delta_2}$ provided that it holds {\em uniformly} in $\Delta_1$ and $\Delta_2$ ({\em i.e.}, provided that $\gamma$ and $\mu$ in Definition \ref{defULB} are independent of $\Delta_1$ and $\Delta_2$).
\end{rmk}

\begin{rmk}\label{rmkcorol}
An interesting corollary of Theorem \ref{theo01} is that the cascade of two USAS systems remains USAS, roughly, provided that the assumptions of Theorem \ref{theo01}, hold with $\delta_1=\delta_2=0$ and that \rref{pilot2} is replaced by the stronger condition:
\begin{equation}\label{newVdot}
\frac{\partial V_1}{\partial t}(t,x_1)+\frac{\partial V_1}{\partial x_1}(t,x_1)f_1(t,x_1,\theta_1) \leq -k_{\Delta_1}V_1(t,x_1)\,,\qquad \forall \norm{x_1}\geq \delta_1\,,
\end{equation}
A similar adaptation of Proposition 2 may be obtained. The proof of this statement is omitted for lack of space, but it follows along the same lines\footnote{The application of Lemma \ref{lemmachangingV} (see Section \ref{prooftheo01}) is no longer required in view of \rref{newVdot}.} as the proof of Theorem \ref{theo01}. Furthermore, it is worth pointing out that, for an \emph{autonomous} $x_1-$subsystem ({\em i.e.}, $f_1(x_1,\theta_1)$), the original requirement \rref{pilot2} remains sufficient.
\end{rmk}

We present a result that may help to check Assumption \ref{A4} in some particular contexts. It requires the non-positivity of the derivative of a Lyapunov function on a sufficiently large domain. The proof is given in Section \ref{sectionproofpropbound}.

\begin{propos}\label{propbound}
Let $b$ be a positive constant. Suppose that there exists a continuously differentiable function $V$ and two class $\cK_\infty$ functions $\underline \alpha$ and $\overline\alpha$ such that, for all $t\in\mRp$ and all $x\in\mathbb R^n$,
\begin{equation}\label{table1}
\underline \alpha(\norm x) \leq V(t,x) \leq \overline \alpha(\norm{x})
\end{equation}
\begin{equation}\label{table2}
x\in\mathcal H(a,b)\quad\Rightarrow\quad \frac{\partial V}{\partial t}(t,x) + \frac{\partial V}{\partial x}(t,x)f(t,x) \leq 0\,,
\end{equation}
where $a$ designates a positive number such that $\overline\alpha(a)<\underline\alpha(b)$. Then, for all $t_0\in\mRp$, the solutions of \rref{sys1} satisfy
$$
\norm{x_0}\leq \overline\alpha^{-1}\circ \underline \alpha(b)\quad\Rightarrow\quad \norm{x(t,t_0,x_0)}\leq b\,,\qquad \forall t\geq t_0\,.
$$
\end{propos}

\section{Application in robot control}\label{sec:eg}

To illustrate the utility of our theorems, we consider the common problem of set-point control of a rigid-joint robot manipulator under PID control and and taking into account the dynamics of the actuators. The Lagrangian dynamics of a robot manipulator with $n$ rigid-joints is given by 
\begin{equation}\label{mod:robot}
  D(q)\ddq + C(q,\dq)\dq + g(q) = u
\end{equation}
where $D(q)\in\mR^{n\times n}$ is symmetric positive definite for all $q\in\mR^n$, $N(q,\dq) := \dot D(q) - 2C(q,\dq)$ is skew-symmetric for all $(q,\,\dq) \in\mR^{n}\times\mR^n$, $u\in\mR^n$ corresponds to the input torques. As most common in the literature of robot control, we restrict our attention to systems satisfying the following ({\em cf.} \cite{SPOVID2,KELLYSBOOK}). 
\begin{standass}\label{standass}
The functions $D(\cdot)$, $C(\cdot,\cdot)$, $g(\cdot)$ are at least twice continuously differentiable and the partial derivatives of their elements are over-bounded by non decreasing functions of\footnote{Notice that this is true for matrices containing polynomial and trigonometric functions which is fairly common in Lagrangian models of physical systems other than mechanical --cf. \cite{THORB,PBCELS}.} $\norm{q}$ and $\norm{\dot q}$. Furthermore, we assume that there exist positive constants $d_m$, $d_M$ and $k_c$ such that\footnote{This is true for instance for open kinematic chains with only revolute or only prismatic joints. See e.g. \cite{SPOVID2,SCISIC}.} for all $q$ and $\dot q\in\mR^n$,
$$
\dty d_m \leq \norm{D(q)} \leq d_M \,,\quad \dty  \norm{C(q,\dot q)}\leq k_c\norm{\dq} \,, \quad \dty \norm{\jac{g(q)}{q}} \leq k_g\,. $$
\end{standass}

We consider that the input torques $u\in\mR^n$ are delivered by Direct-Current (DC) motors, whose dynamics are given by 
\begin{equation}\label{eq:motor}
  L \dot i + R i + k_b \dot q = v
\end{equation}
where  $i\in\mR^n$ is the vector of rotor currents, $L$ and $R$ are the rotors' inductances and resistances respectively\footnote{For simplicity and without loss of generality we consider that we have $n$ {\em identical} motors --{\em i.e.} same resistance, inductance, torque constant, etc.},  $k_b\dot q_j$ with $j\leq n$ is the back-emf voltage in each motor and $v$ is the vector of input voltages, {\em i.e.}, the control inputs. Each motor produces an output torque $u_j = k_t i_j$ with $k_t>0$. 

Our control problem is to design $v$ so that the robot manipulator stabilizes at a desired constant set-point $(q=q^*$, $\dot q =0)$. The aim is that the robot coordinates approach the reference operating point from any initial conditions in an \emph{arbitrarily large} set (limited only by physical constraints). Furthermore, it is imposed that control be of the PID type. Accordingly, disregarding the DC motor dynamics, the input torques that achieve the control objective are given by 
\begin{subequations}\label{eq:PID}
\begin{eqnarray}\label{eq:PID:a}
  u^* & = & -k_p\tq - k_d \dq + \nu\\
\label{eq:PID:b}
  \dot \nu & = & -k_i \tq\,, \qquad \nu(0) := \hat g(q^*)
\end{eqnarray}
\end{subequations}
where $k_p$, $k_d$ and $k_i$ are positive design control gains  and $\hat g(q^*)$ is the best ``guess''  of the {\em unknown} constant pre-computed gravitational forces vector. 

We stress that the above setting is fairly common in practice of robot control; not only PID control is probably the most popular control technique but, often, industrial manipulators come with a black-box controller of PID type, meaning that control design for the user of an industrial robot boils down to gain-tuning for the built-in PID. 

Here, our control objective is achieved via cascaded-based control, relying on our results on USPAS. The approach consists in designing a reference $i^*:= u^*/k_t$ (so that, when $\tilde i := i - i^* = 0$, we have that $u=u^*$) and building a control law $v$ that makes $\tilde i$ go to zero. 

\begin{proposition}\label{prop:USPASofrobot}
Under Standing assumption \ref{standass}, the system \rref{mod:robot}, \rref{eq:motor} in closed-loop with \rref{eq:PID} and 
$$
v  :=  R' \tilde i + Ri^* + k_b \dot q  + L\dot i^*\,,\quad i^*=\frac{u^*}{k_t}\,.
$$
is uniformly semiglobally asymptotically stable (USAS). 
\end{proposition}
Proposition \ref{prop:USPASofrobot} establishes that, if one knows how to semiglobally asymptotically stabilize a robot using PID control when neglecting the DC drive dynamics, then the same stability property can be established when these dynamics are taken into account. In other words, we claim that, given any domain of initial errors, one can always find control gains (namely $k_p$, $k_i$, $k_d$) such that the point $(q=q^*,\dot q=0)$ is uniformly asymptotically stable on this set of initial conditions. Moreover, the tuning of the PID control gains can be made disregarding the DC drive dynamics. 

\begin{sketchproof}\emph{(Proposition \ref{prop:USPASofrobot})}
For analytical purposes, let $\ep_1>0$ be sufficiently small and define the variables $s:= \frac{1}{\ep_1}\tq + \frac{1}{k_i}(g(q^*)-\nu)$ and $k_p' := k_p - \frac{k_i}{\ep_1} > 0$. Then, the closed-loop system can be written in the following cascaded form:
\begin{subequations}\label{robot:cascade}
\begin{eqnarray}\label{robot:cascade:a}
 D(q)\ddq + C(q,\dq)\dq + g(q) - g(q^*) + k_p'\tq + k_d \dq - k_i s & = & k_t \tilde i \\
\label{robot:cascade:b}
 \dot {s} & = & \tq + \frac{1}{\ep_1} \dot q
\end{eqnarray}
\end{subequations}
\begin{equation}
\shiftright{87mm}\dot{\tilde i} \ \ =\ \ -\frac{R+R'}{L}\tilde i \,.\label{motor:cascade}
\end{equation}

Then, the proof of the proposition can be constructed by applying Theorem \ref{theo01} and using Remark \ref{rmkcorol}. For this, three basic properties must be shown: 1) the motor closed-loop system \rref{motor:cascade} is USAS, 2) the robot system in closed loop with the PID controller, {\em i.e.} \rref{robot:cascade} with $\tilde i \equiv 0$, is USAS and 3) the trajectories of the cascaded system are uniformly bounded. 

\noindent \underline{Sketch of proof of 1)}: This follows directly by observing that \rref{motor:cascade} is uniformly {\em globally} exponentially stable for any positive $R'$.

\noindent \underline{Sketch of proof of 2)}:  This follows after lengthy calculations that we do not include here for lack of space and because they follow along similar prooflines as those in \cite{PI2D}. For the purpose of illustrating the use of Theorem \ref{theo01} we claim that the result can be established using the Lyapunov function 
$$
V_1 := \frac{1}{2}\dot q^\top D(q) \dot q + \frac{1}{2}k_p'\norm{\tq}^2 + U(q) - U(q^*) - \tilde q^\top g(q^*) + \frac{\ep_1 k_i}{2}s^2 + \ep_1 \tq^\top D(q) \dot q + \ep_2 s^\top  D(q) \dot q
$$
where $U(\cdot)$ corresponds to the gravitational potential energy function and $\varepsilon_2>0$. We emphasize that, except for the last term and a slightly different notation, the function above comes from \cite{PI2D} and, more indirectly, from energy-like Lyapunov functions widely used in the robot-control literature. In particular, conditions for positive definiteness\footnote{Positive definiteness of $V$ may be established if $k_p'\geq k_g$ and if $\ep_1$, $\ep_2$ are sufficiently small.} and radial unboundedness hold under the Standing assumption 1 and can be derived following \cite{PI2D} and some of the references therein. 

It can also be shown that, for any $\Delta>0$, there exist $k_p'$, $k_d$ and $k_i$ such that the time derivative of $V_1$ along the trajectories of  \rref{robot:cascade} satisfies, for all $\norm{x_1}:=\sqrt{\norm{\dq}^2 + \norm{\tq}^2  + \norm{s}^2} \leq \Delta_1$ and all $t\geq 0$,
$$
\dot V \leq -\frac{k_d}{2}\norm{\dot q}^2 - \frac{\ep_1 k_p'}{2}\norm{\tq}^2 - \frac{\ep_2 k_i}{2}s^2\,,
$$
provided that $\ep_1$ and $\ep_2$ are chosen small enough, and that the choice of the gains can be made in the following manner\footnote{We can actually show that the $\mathcal D$-set ({\em cf.} Definition \ref{Dset}) of \rref{robot:cascade} is $$\mathcal D(0,\Delta_1)=\left\{(k_d,k_p,k_i)\in\mR^3~:~k_d\geq a_d+b_d\Delta_1,\,k_p'\geq a_p+b_p\Delta_1,\,k_i\geq a_i+b_i\Delta_1\right\}\,.$$}: 
\begin{equation}\label{gains}
k_d=a_d+b_d\Delta_1\,,\quad k_p'=a_p+b_p\Delta_1\,,\quad k_i=a_i+b_i\Delta_1\,,
\end{equation}
where $a_d$, $a_p$, $a_i$, $b_d$, $b_p$, $b_i$ are positive constants. In addition, the bound on the gradient \rref{pilot3} directly follows from the smoothness of $V_1$ and its time-independency. This establishes Assumption \ref{A1} when considering $\delta_1=0$.

\noindent \underline{Sketch of proof of 3)}: This can be established following the same analysis as above, with the same Lyapunov function, but considering the interconnection term $k_t \tilde i$ on the right hand side of the closed loop equation \rref{robot:cascade:a}. One obtains that, in this case,  
$$
\dot V \leq -\frac{k_d}{2}\norm{\dot q}^2 - \frac{\ep_1 k_p'}{2}\norm{\tq}^2 - \frac{\ep_2 k_i}{2}s^2 +
 \left(\norm{\dq}+ \ep_1\norm{\tq} + \ep_2\norm{s}\right)k_t \norm{\tilde i}\,.
$$
which, in view of the uniform boundedness of $\tilde i(t)$, satisfies $\dot V \leq 0$ for ``large'' values of $\norm{x_1}$. This observation, combined with the linear dependency of the gains \rref{gains} in $\Delta_1$, allows to fulfill Assumption \ref{A4}, in view of Proposition \ref{propbound}.
\end{sketchproof}
\begin{remark} 
Cascaded-based control of robots taking account of the robot dynamics was first used in \cite{PANORT} in trajectory control of manipulators with AC drives. In that problem, the motor dynamics is highly nonlinear and global asymptotic stability of the closed-loop system of \rref{mod:robot} with the corresponding ideal control input $u^*$ is obtained. However, for the problem that we address here, we are not aware of any proof of global asymptotic stability of the closed-loop \rref{mod:robot} with PID control hence global results for cascaded systems fail in this setting. See also \cite{AILOGI97} for a result on control of robots taking into account the DC motors' dynamics, but with knowledge of the gravity terms $g(q)$.
\end{remark}


\section{Proofs}
\label{sec:proofs}

\subsection{Proof of Theorems \ref{theo01}}\label{prooftheo01}
\vspace{-0.5cm}We start by introducing the following result, which is a direct adaptation of \cite[Proposition 13]{PRAWAN} and allows $V_1$ to be transformed into a more convenient form. 
\begin{lemma}\label{lemmachangingV}
Let $\delta$ be a nonnegative constant and $X$ be a subset of $\mR^n\setminus \mathcal B_\delta$. Suppose that there exist a continuously differentiable function $V:\mRp\times X\to\mRp$ and some class $\cK_\infty$ functions $\underline\alpha$, $\overline\alpha$, $\alpha$ such that, for all $x\in X$ and all $t\geq 0$, 
$$
\underline\alpha(\norm {x}_{\delta}) \leq V(t,x) \leq \overline\alpha(\norm {x})
$$
$$
\frac{\partial V}{\partial t}(t,x)+\frac{\partial V}{\partial x}(t,x)f(t,x) \leq -\alpha(\norm x)\,.
$$
Then, for any positive $k$, there exists a continuously differentiable function $\mathcal V:\mRp\times X\to\mRp$ and class $\cK_\infty$ functions $\tilde{\underline\alpha}$, $\tilde{\overline\alpha}$ such that, for all $x\in X$ and all $t\geq 0$, 
\begin{equation}\label{ventilo1}
\tilde{\underline\alpha}(\norm {x}_{\delta}) \leq \mathcal V(t,x) \leq \tilde{\overline\alpha}(\norm {x})
\end{equation}
\begin{equation}\label{ventilo2}
\frac{\partial \mathcal V}{\partial t}+\frac{\partial \mathcal V}{\partial x}f(t,x) \leq -k\mathcal V\,,
\end{equation}
and, for any $s\in\mRp$, it holds that 
\begin{equation}\label{ventilo2bis}
\tilde{\underline\alpha}^{-1}\circ\tilde{\overline\alpha}(s)=\underline\alpha^{-1}\circ\overline\alpha(s)\,.
\end{equation}
If, in addition, there exists a continuous nondecreasing function $c:\mRp\to\mRp$ such that, for all $x\in X$ and all $t\geq 0$,
$$
\norm{\frac{\partial V}{\partial x}(t,x)}\leq c(\norm{x})\,,
$$
then there exists a continuous nondecreasing function $\tilde c:\mRp\to\mRp$ such that, for all $x\in X$ and all $t\geq 0$, 
\begin{equation}\label{ventilo3}
\norm{\frac{\partial \mathcal V}{\partial x}(t,x)}\leq \tilde c(\norm{x})\,.
\end{equation}
\end{lemma}

\begin{proof}
Following the prooflines of \cite[Proposition 13]{PRAWAN}, we see that the function $\mathcal V$ can be defined as $\rho\circ V$ where 
$$
\left\{\begin{array}{ll} \rho(s)&=\,\exp\left(\int_1^s \frac{2dq}{a(q)}\right)\,,\quad \forall s>0 \\ \rho(0)&=\,0\,, \end{array}\right.
$$
and $a$ is a convenient class $\cK$ function. The bound \rref{ventilo2} can be established following the same reasoning as in the proof of \cite[Proposition 13]{PRAWAN}. Furthermore, \rref{ventilo1} can be satisfied with $\tilde{\underline\alpha}:=\rho\circ \underline\alpha$ and $\tilde{\overline\alpha}:=\rho\circ \overline\alpha$ as $\rho\in\cKinfty$, and we therefore have that
$$
\tilde{\underline\alpha}^{-1}\circ\tilde{\overline\alpha}(s)=\left(\rho\circ\underline\alpha\right)^{-1}\circ\left(\rho\circ\overline\alpha\right)(s) = \left(\underline\alpha^{-1}\circ\rho^{-1}\right)\circ\left(\rho\circ\overline\alpha\right)(s)= \underline\alpha^{-1}\circ\overline\alpha(s)\,.
$$ 
Concerning the bound on the gradient we have that, for all $x\in X$ and all $t\in\mRp$,
$$
\norm{\frac{\partial \mathcal V}{\partial x}(t,x)}\leq \frac{2\mathcal V(x)}{a(V(x))}\norm{\frac{\partial V}{\partial x}(t,x)} \leq \frac{2\tilde{\overline\alpha}(\norm x)}{a(\underline\alpha(\norm x))}c(\norm{x}) \leq \tilde c(\norm x)\,,
$$
where $\tilde c(s):= \frac{2\tilde{\overline\alpha}(s)}{a(\underline\alpha(\delta))}c(s)$\ac{A counter-example for a function $a$ satisfying $a\in\cK\cap C^1$, $a(s)\leq s$, $a'(0)$ and $$\lim_{s\to 0}\frac{1}{a(s)}\exp\left(\int_1^s\frac{d\tau}{a(\tau)}\right)<\infty$$ is $a(s)=\ln\left(\norm{\ln(s)}\right)$.}, which establishes the result.
\end{proof}

Consider the function $V_1$ generated by Assumption \ref{A1} and let Lemma \ref{lemmachangingV} with $x=\mathcal H(\Delta_1,\delta_1)$ generate a function $\mathcal V_1$, class $\cK_\infty$ functions $\tilde{\underline\alpha}_{\delta_1,\Delta_1}$, $\tilde{\overline\alpha}_{\delta_1,\Delta_1}$, a positive constant $k_{\delta_1,\Delta_1}$  and a continuous nondecreasing function $\tilde c_{\delta_1,\Delta_1}$ such that, for all $x\in \mathcal H(\delta_1,\Delta_1)$ and all $t\in\mRp$, \rref{ventilo1}-\rref{ventilo3} hold for $\mathcal V_1$. From \rref{condadd3} and \rref{ventilo2bis}, it also holds that, for any $\Delta_1>0$,
\begin{equation}\label{condadd4}
\lim_{\delta_1\to 0}\tilde{\underline\alpha}_{\delta_1,\Delta_1}^{-1}\circ\tilde{\overline\alpha}_{\delta_1,\Delta_1}(\delta_1)=0\,.
\end{equation}
In the sequel, in order to lighten the notations, we refer to $\tilde{\underline\alpha}_{\delta_1,\Delta_1}$ as simply $\underline \alpha_1$, $\tilde{\overline\alpha}_{\delta_1,\Delta_1}$ as $\overline \alpha_1$, $k_{\delta_1,\Delta_1}$ as $k_1$ and $\tilde c_{\delta_1,\Delta_1}$ as $c_1$. Even though no longer explicit with these new notations, the dependency of these functions in $\delta_1$ and $\Delta_1$ should be kept in mind. For any given positive $\delta_1$, $\Delta_1$, $\delta_2$ and $\Delta_2$ satisfying $\Delta_1>\max\{\delta_1,\Delta_0\}$ and $\Delta_2>\delta_2$, let $\gamma_{\Delta_1,\Delta_2}$ be generated by Assumption \ref{A4} and define
\begin{equation}\label{Delta}
\Delta:=\min\left\{\Delta_1\,;\,\Delta_2\,;\,\gamma_{\Delta_1,\Delta_2}(\Delta_1,\Delta_2)\right\}\,.
\end{equation}
Next, choose any $\theta_1^\star\in\Theta_1$ satisfying Assumption \ref{A1} and any $\theta_2^\star\in\mathcal D_{f_2}(\delta_2,\Delta_2)\cap\Theta_2$ given by Assumption \ref{A4}. We show that, provided that $\delta_1$, $\delta_2$ are  sufficiently small and that $\Delta_1$, $\Delta_2$ are large enough, there exists $\delta\in (0;\Delta)$ such that  $\mathcal B_\delta$ is UAS on $\mathcal B_\Delta$ for the system \rref{syscorol} with $\theta^\star=(\theta_1^\star,\theta_2^\star)$. To that end, we first show that there exists  a positive $\delta_3$ such that the ball $\mathcal B_{\delta_3}$ is uniformly stable. More precisely, we construct $\eta\in \cKinfty$ and $\delta_3>0$  such that, for all $x_0\in\mathcal B_\Delta$, 
\begin{equation}\label{boundx14bis}
\norm{x_1(t,t_0,x_{0},\theta^\star)}_{\delta_3}\leq \eta(\norm{x_0})\,.
\end{equation}
Then, we use this property to prove that a ball, larger than $\mathcal B_{\delta_3}$, is UA on $\mathcal B_\Delta$ and we construct a $\cKL$ estimate for the solutions. Finally, we show that the estimates of the domain of attraction and of the ball to which solutions converge can be arbitrarily enlarged and diminished respectively.

\subsubsection{Proof of uniform stability of a ball}
\vspace{-0.5cm}The total time derivative of ${\mathcal V}$ along the trajectories of \rref{syscorol} with $\theta=\theta^\star$ yields
\begin{displaymath}
\dot {\mathcal V}_1=\frac{\partial {\mathcal V}_1}{\partial t}+\frac{\partial {\mathcal V}_1}{\partial x_1}\big(f_1(t,x_1,\theta_1^\star)+g(t,x,\theta^\star)x_2\big).
\end{displaymath}
Therefore it holds that, for all $x_1\in\mathcal H(\delta_1,\Delta_1)$ and all $t\geq 0$
\begin{eqnarray}
\dot {\mathcal V}_1&\leq& -k_1{\mathcal V}_1 + \norm{\frac{\partial {\mathcal V}_1}{\partial x_1}}\norm{g(t,x,\theta^\star)}\norm{x_2}\nonumber \\&\leq& -k_1{\mathcal V}_1 + c_1(\norm{x_1})G(\norm{x})\norm{x_2}\,.\nonumber
\end{eqnarray}
Defining 
$$
\Gamma:=\left\{ t\geq  t_0~|~\delta_1\leq\norm{x_1(t,t_0,x_{0},\theta^\star)}\leq \Delta_1\right\}\,,
$$
and using the shorthand notation $x_1(t)$  for $x_1(t,t_0,x_{0},\theta^\star)$ and $v_1(t):={\mathcal V}_1(t,x_1(t))$ we get that, for any $x_0\in\mathcal B_\Delta$ and any $t\in\Gamma$, 
$$
\dot v_1(t)\leq -k_1v_1(t) + c_1(\norm{x_1(t)})G(\norm{x(t)})\norm{x_2(t)}\,.
$$
From Assumption \ref{A4}, and in view of \rref{Delta}, we can see that, for all $x_0\in\mathcal B_\Delta$, it holds that $x(t)\in\mathcal B_{\Delta_1}$. Hence, for all $t\in\Gamma$,
$$
\dot v_1(t)\leq -k_1v_1(t) + c_1(\Delta_1)G(\Delta_1)\norm{x_2(t)}\,.
$$
Let Assumption \ref{A2} generate a class $\mathcal{KL}$ function $\beta_2$ such that\footnote{Notice that $\norm{s}_a\leq b \Leftrightarrow \norm{s}\leq a+b$, $\forall s\in\mR^n$, $a, b\geq0$.} for any $x_{20}\in\mathcal B_{\Delta_2}$ and any $t\geq t_0$,
$$
\norm{x_2(t)}\leq \beta_2(\norm{x_{20}},t-t_0)+\delta_2\,.
$$
It follows that, for all $x_0\in \mathcal B_{\Delta}$ and all $t\in\Gamma$,
\begin{equation}
 \label{V1dot4bis}
\dot v_1(t)\leq -k_1v_1(t) + c_1(\Delta_1)G(\Delta_1)\big[\beta_2(\norm{x_{20}},t-t_0) + \delta_2\big]
\end{equation}
which implies that
\begin{equation} \label{V1dot6}
x(t)\in \mathcal H(\delta_1,\Delta_1)\quad\Rightarrow\quad \dot v_1(t)\leq -k_1v_1(t) + c_3(\norm{x_0})\,,
\end{equation}
with 
$$
c_3(s):=c_1(\Delta_1)G(\Delta_1)(\beta_2(s,0)+\delta_2)\,,\qquad \forall s\geq 0\,.
$$
The rest of the proof of uniform stability consists in integrating \rref{V1dot6} over $\Gamma$, in order to construct a bound like \rref{boundx14bis}. To that end, we introduce the following tool, which may be viewed as a comparison theorem for differential inequalities that hold only out of a ball centered at zero. 
\begin{lemma}\label{lemmaintegrate}
Let $\delta$ be a nonnegative constant and $X$ be a subset of $\mR^n$ containing $\mathcal B_\delta$. Assume that there exists a continuously differentiable function $V:\mRp\times\mR^n\to\mRp$, class $\cK_\infty$ functions $\underline\alpha$ and $\overline\alpha$, a positive constant $k$ and nonnegative constant $c$ such that, for all $x\in X$ and all $t\in \mRp$,
$$
\underline\alpha(\norm{x}_\delta)\leq V(t,x) \leq \overline\alpha(\norm{x})
$$
and, for all $x_0\in\mR^n$ and all $t_0\in \mRp$,
$$
x(t,t_0,x_0)\in X\setminus \mathcal B_\delta \quad\Rightarrow\quad \dot V(t,x(t,t_0,x_0))\leq -k V(t,x(t,t_0,x_0))+c\,.
$$
Then, for all $x_0\in\mR^n$ and $t_0\in\mRp$ such that $x(t,t_0,x_0)\in X$ $\forall t\geq t_0$, we have that
$$
\norm{x(t,t_0,x_0)}_\delta\leq \underline\alpha^{-1}\left(\overline\alpha(\delta)+\frac{c}{k}\right) +\underline\alpha^{-1}\left(\overline\alpha(\norm{x_0})e^{-k(t-t_0)}+\frac{c}{k}\right)\,,\qquad\forall t\geq t_0\,.
$$
\end{lemma}
\begin{proof}
For simplicity, we write $x(\cdot,t_0,x_0)$ as $x(\cdot)$ and we define $v(\cdot):=V(\cdot,x(\cdot))$. We distinguish two cases: whether the trajectories start from outside or inside $\mathcal B_{\delta}$.

\noindent \underline{\sl Case 1}: $\norm{x_0}>\delta$.\\
\noindent In this case, there exists\footnote{If $\norm{x(t)}>\delta$ forever after, we consider that $T_0=\infty$.} $T_0\in(0;\infty]$ such that $\norm{x(t)}>\delta$ for all $[t_0;t_0+T_0)$ and $\norm{x(t_0+T_0)}=\delta$. Hence, using the comparison lemma, we get that $v(t)\leq (v(t_0)-\frac{c}{k})e^{-k(t-t_0)}+\frac{c}{k}$ for all $t\in [t_0;t_0+T_0)$. Using the bounds on $V$, it follows that
$$
\norm{x(t)}_\delta\leq \underline\alpha^{-1}\left(\overline\alpha(\norm{x_0})e^{-k(t-t_0)}+\frac{c}{k}\right)\,,\qquad \forall t\in [t_0;t_0+T_0)\,.
$$
In addition, for each $t\geq t_0+T_0$, either $\norm{x(t)}\leq\delta$ in which case\footnote{\label{fn1}This is direct by noticing that $\underline\alpha(s)\leq\overline\alpha(s)$ for all $s\in\mRp$ and that $c/k\geq 0$.} $\norm{x(t)}_\delta\leq\underline\alpha^{-1}\left( \overline\alpha(\delta)+c/k\right)$, or $\norm{x(t)}>\delta$. In this second case, we can again invoke the continuity of the solution to see that there exists a nonempty time-interval $[\tau;\tau+T]$, with $T\in(0;\infty]$, containing $t$ and such that $\norm{x(s)}>\delta$ for all $s\in(\tau;\tau+T]$, with $\norm{x(\tau)}=\delta$. Hence, integrating from $\tau$ to $t\in[\tau;\tau+T]$, we obtain in the same way as before that, whenever $\norm{x(t)}>\delta$, it holds that
\begin{equation}
\norm{x(t)}_\delta \leq \underline \alpha^{-1}\left(\overline \alpha(\delta)e^{-k(t-\tau)}+\frac{c}{k}\right)
\leq \underline \alpha^{-1}\left(\overline \alpha(\norm{x_0})e^{-k(t-t_0)}+\frac{c}{k}\right)\,.\label{eq2045}
\end{equation}
To sum up, for all $t\geq t_0$, we have the following:
\begin{equation}\label{eq2046}
\norm{x_0}>\delta\quad\Rightarrow\quad \norm{x(t)}_\delta \leq \underline \alpha^{-1}\left(\overline \alpha(\norm{x_0})e^{-k(t-t_0)}+\frac{c}{k}\right)\,.
\end{equation}

\noindent \underline{\sl Case 2}: $\norm{x_0}\leq\delta$.\\
\noindent In this case, as long as $\norm{x(t)}\leq \delta$, we trivially$\hspace{-2mm}~^{\ref{fn1}}$ have that $\norm{x(t)}_\delta\leq \underline\alpha^{-1}\left(\overline\alpha(\delta)+c/k\right)$. If $\norm{x(t)}>\delta$ at some instant $t>t_0$, then, again, there exists a nonempty time-interval $[\tau;\tau+T]$, with $T\in(0;\infty]$ and $\tau>t_0$, containing $t$ and such that $\norm{x(s)}>\delta$ for all $s\in(\tau;\tau+T]$, with $\norm{x(\tau)}=\delta$. Thus, from \rref{eq2045}, we obtain that
$$
\norm{x(t)}_\delta \leq \underline \alpha^{-1}\left(\overline \alpha(\delta)e^{-k(t-\tau)}+\frac{c}{k}\right)\leq \underline \alpha^{-1}\left(\overline \alpha(\delta)+\frac{c}{k}\right)\,.
$$
Hence, for all $t\geq t_0$,
\begin{equation}\label{eq2047}
\norm{x_0}\leq\delta\quad\Rightarrow\quad\norm{x(t)}_{\delta} \leq \underline \alpha^{-1}\left(\overline \alpha(\delta)+\frac{c}{k}\right)\,.
\end{equation}
The conclusion follows from \rref{eq2046} and \rref{eq2047}.
\end{proof}

Applying Lemma \ref{lemmaintegrate} to \rref{V1dot6} with $V=V_1$, $k=k_1$, $c=c_3(\norm{x_0})$ and $X=\mathcal B_\Delta$, we get in view of Assumption \ref{A4} and \rref{Delta} that, for all $x_0\in\mathcal B_\Delta$ and all $t_0\in\mRp$,
$$
\norm{x(t)}_{\delta_1}\leq \underline\alpha_1^{-1}\left(\overline\alpha_1(\delta_1)+\frac{c_3(\norm{x_0})}{k_1}\right) +\underline\alpha_1^{-1}\left(\overline\alpha_1(\norm{x_0})+\frac{c_3(\norm{x_0})}{k_1}\right)\,,\qquad\forall t\geq t_0\,.
$$
Define the following:
\begin{eqnarray}
\delta_3 &:=& \delta_1+\underline \alpha_1^{-1}\left(\overline \alpha_1(\delta_1)+\frac{c_3(0)}{k_1}\right)+ \underline \alpha_1^{-1}\left(\frac{c_3(0)}{k_1}\right) \nonumber \\
&=& \delta_1+\underline \alpha_1^{-1}\left(\overline \alpha_1(\delta_1)+\frac{c_1(\Delta_1)G(\Delta_1)\delta_2}{k_1}\right)+ \underline \alpha_1^{-1}\left(\frac{c_1(\Delta_1)G(\Delta_1)\delta_2}{k_1}\right)\,.\label{delta3}
\end{eqnarray}
and, for all $s\in\mRp$,
$$
\eta(s):= \underline \alpha_1^{-1}\left(\overline \alpha_1(\delta_1)+ \frac{c_3(s)}{k_1}\right)+ \underline \alpha_1^{-1}\left(\overline \alpha_1(s)+ \frac{c_3(s)}{k_1}\right)- \underline \alpha_1^{-1}\left(\overline \alpha_1(\delta_1)+\frac{c_3(0)}{k_1}\right) -\underline \alpha_1^{-1}\left(\frac{c_3(0)}{k_1}\right)\,.
$$
We then conclude that, for any $x_0\in\mathcal B_\Delta$ and all $t_0\in\mRp$, it holds that
\begin{equation} \label{boundx12}
\norm{x_1(t)}_{\delta_3}\leq \eta(\norm{x_0})\,,\qquad\forall t\geq t_0\,.
\end{equation}
Uniform stability of $\mathcal B_{\delta_3}$ on $\mathcal B_{\Delta}$ follows by noticing that $\eta$ is a class $\cK$ function. This can be seen by recalling that $c_3$ is a continuous increasing function.

\subsubsection{Proof of uniform attractiveness of a ball}
\vspace{-0.5cm}Consider again \rref{V1dot4bis}. Since $\beta_2$ is a $\mathcal{KL}$ function there is a time $t_1\geq 0$, independent of $t_0$ and $x_0$, such that
\begin{displaymath}
\beta_2(\Delta,t-t_0)\leq \delta_2\,,\qquad \forall\, t\geq t_0+t_1\,.
\end{displaymath}
Hence \rref{V1dot4bis} implies that, for all $t\in\Gamma\cap\mR_{\geq t_0+t_1}$ and all $x_0\in\mathcal B_\Delta$,
$$
\dot v_1(t)\leq -k_1 v_1(t) + 2c_1(\Delta_1)G(\Delta_1)\delta_2\,.
$$
Applying again Lemma \ref{lemmaintegrate} and recalling that, from Assumption \ref{A4}, $\norm{x_1(t_0+t_1)}\leq \Delta_1$, it follows that, for all $x_0\in\mathcal B_\Delta$, all $t_0\in\mRp$ and all $t\geq t_0+t_1$,
\begin{eqnarray*}
\norm{x(t)}_{\delta_1}&\leq& \underline\alpha_1^{-1}\left(\overline\alpha_1(\delta_1)+\frac{2c_1(\Delta_1)G(\Delta_1)\delta_2}{k_1}\right) +\underline\alpha_1^{-1}\left(\overline\alpha_1(\norm{x(t_0+t_1)})e^{-k_1(t-t_0-t_1)}+\frac{2c_1(\Delta_1)G(\Delta_1)\delta_2}{k_1}\right)\\
&\leq& \underline\alpha_1^{-1}\left(\overline\alpha_1(\delta_1)+\frac{2c_1(\Delta_1)G(\Delta_1)\delta_2}{k_1}\right) +\underline\alpha_1^{-1}\left(\overline\alpha_1(\Delta_1)e^{-k_1(t-t_0-t_1)}+\frac{2c_1(\Delta_1)G(\Delta_1)\delta_2}{k_1}\right)\,.
\end{eqnarray*}
Defining
$$
t_2:=t_1+\frac{1}{k_1}\ln\left(\frac{\overline \alpha_1(\Delta_1)}{\overline \alpha_1(\delta_1)}\right)\,,
$$
we then see that, for all $x_0\in\mathcal B_\Delta$,
\begin{equation}\label{delta4}
\norm{x_1(t)} \leq \delta_4:=\delta_1+2\underline \alpha_1^{-1}\left(\overline \alpha_1(\delta_1) +\frac{2c_1(\Delta_1)G(\Delta_1)\delta_2}{k_1}\right)\,,\qquad \forall t\geq t_0+t_2\,. 
\end{equation}
In other words, we have that 
$$
\norm{x_1(t,t_0,x_{10})}_{\delta_4} = 0\,,\qquad \forall t\geq t_0+t_2\,.
$$ 
Finally, let 
\begin{equation}\label{delta}
\delta:=\max\left\{\delta_2\,;\,\delta_3\,;\,\delta_4\right\}\,,
\end{equation}
Then we see that \rref{boundx12} implies that $\norm{x_1(t)}_{\delta}\leq \eta(\norm{x_0})$ for all $t\geq t_0$. From this and what precedes it is not hard to see that, for all $x_0\in\mathcal B_\Delta$,
$$
\norm{x_1(t)}_{\delta}\leq \eta(\norm{x_0})e^{-(t-t_0-t_2)}\,,\qquad \forall t\geq t_0\,.
$$
Thus, recalling that $t_2$ depends neither on $t_0$ nor on $x_0$, and defining
$$
\beta(s,t):=\max\left\{ \eta(s)e^{-(t-t_2)}\,;\,\beta_2(s,t)\right\}\,,\qquad \forall s,t\geq 0\,,
$$
we conclude that, for all $x_0\in\mathcal B_\Delta$,
$$
\norm{x(t)}_\delta\leq \beta(\norm{x_0},t-t_0)\,,\qquad\forall t\geq t_0\,.
$$
UAS of $\mathcal B_\delta$ on $\mathcal B_\Delta$ follows by noticing that $\beta$ is a class $\cKL$ function.

\subsubsection{Semiglobality and practicality}
\vspace{-0.5cm}It is only left to show that $\delta$ and $\Delta$ can be arbitrarily reduced and enlarged respectively. It follows directly from \rref{condgamma} and $\rref{Delta}$ that, by picking $\Delta_1$ and $\Delta_2$ large enough, $\Delta$ can be made arbitrarily large. 

Concerning $\delta$, we see with \rref{delta3} and \rref{delta4} that
$$
\delta_3\leq\delta_4=\delta_1+2\underline \alpha_1^{-1}\left(\overline \alpha_1(\delta_1) +\frac{2c_1(\Delta_1)G(\Delta_1)\delta_2}{k_1}\right)\,.
$$
Hence, in view of \rref{condadd4} and recalling that $c_1$, $G$, $k_1$ and $\underline \alpha_1$ are independent of $\delta_2$, we see that, for the chosen $\Delta_1$ and $\Delta_2$, both $\delta_3$ and $\delta_4$ can be taken as small as wanted by picking $\delta_1$ and $\delta_2$ sufficiently small. Hence, in view of \rref{delta}, it is also the case for $\delta$.

Thus, it suffices to pick the parameters $\theta_1^\star$ and $\theta_2^\star$ generated by the chosen $\delta_1$, $\Delta_1$, $\delta_2$ and $\Delta_2$, to conclude that, for any $\Delta>\delta>0$, there exists some parameters $\theta_1^\star\in\Theta_1$ and $\theta_2^\star\in\Theta_2$ such that $\mathcal B_\delta$ is UAS on $\mathcal B_\Delta$ for the system \rref{syscorol} with $\theta=\theta^\star$, which establishes the result.


\subsection{Proof of Proposition \ref{prop1}} \label{proofprop1}

\vspace{-0.5cm} 
We prove this result by showing that it actually constitutes a special case of Theorem \ref{theo01}. To this end, first notice that applying Lemma \ref{lemmachangingV} (see Section \ref{prooftheo01}) to $V$ ensures the existence of a continuously differentiable function $\mathcal V$ such that, for all $x\in\mathcal H(\delta,\Delta)$ and all $t\in\mRp$,
$$
\tilde{\underline\alpha}_{\delta,\Delta}(\norm {x}_{\delta}) \leq \mathcal V(t,x) \leq \tilde{\overline\alpha}_{\delta,\Delta}(\norm {x})
$$
\begin{equation}\label{ventilo2bis}
\frac{\partial \mathcal V}{\partial t}(t,x)+\frac{\partial \mathcal V}{\partial x}(t,x)f(t,x,\theta) \leq -k_{\delta,\Delta}\mathcal V(t,x)
\end{equation}
$$
\norm{\frac{\partial \mathcal V}{\partial x}(t,x)}\leq \tilde c_{\delta,\Delta}(\norm{x})
$$
hold with a positive $k_{\delta,\Delta}$, a continuous nondecreasing function $\tilde c_{\delta,\Delta}$ and some class $\cK_\infty$ functions $\tilde{\underline \alpha}_{\delta,\Delta}$ and $\tilde{\overline \alpha}_{\delta,\Delta}$ satisfying
$$
\tilde{\underline \alpha}^{-1}_{\delta,\Delta} \circ \tilde{\overline \alpha}_{\delta,\Delta}(s) = \underline \alpha^{-1}_{\delta,\Delta} \circ \overline \alpha_{\delta,\Delta}(s)\,,\qquad \forall s\geq 0\,.
$$
Inverting the two sides of this inequality yields:
$$
\tilde{\overline \alpha}^{-1}_{\delta,\Delta} \circ \tilde{\underline \alpha}_{\delta,\Delta}(s) = \overline \alpha^{-1}_{\delta,\Delta} \circ \underline \alpha_{\delta,\Delta}(s)\,,\qquad \forall s\geq 0\,.
$$
Consequently, in view of \rref{condadd} and \rref{condadd2}, we have that, for all $\Delta>0$,
$$
\lim_{\delta\to 0} \tilde{\underline \alpha}^{-1}_{\delta,\Delta} \circ \tilde{\overline \alpha}_{\delta,\Delta}(\delta) = 0
$$
and, for all $\delta>0$,
\begin{equation}\label{condalpha}
\lim_{\Delta\to \infty}\tilde{\overline \alpha}^{-1}_{\delta,\Delta} \circ \tilde{\underline \alpha}_{\delta,\Delta}(\Delta)=\infty\,.
\end{equation}
Based on this, let $\Delta$ be any given positive constant and choose $\delta$ small enough that
$$
\tilde{\underline \alpha}^{-1}_{\delta,\Delta} \circ \tilde{\overline \alpha}_{\delta,\Delta}(\delta)<\Delta\,.
$$
Then, the requirements of Proposition \ref{propbound} are fulfilled and we get that
$$
\norm{x_0}\leq \tilde{\overline \alpha}^{-1}_{\delta,\Delta} \circ \tilde{\underline \alpha}_{\delta,\Delta}(\Delta)\quad\Rightarrow\quad \norm{x(t)}\leq \Delta\,,\qquad \forall t\geq t_0\,.
$$
In view of \rref{condalpha}, we then see that the solutions of \rref{sys2} satisfy a uniform boundedness as in Assumption \ref{A4}. Moreover, from \rref{ventilo2bis}, and defining $v(t):=\mathcal V(t,x(t,t_0,x_0,\theta))$ and $\Gamma:=\left\{t\geq t_0~|~x(t,t_0,x_0,\theta)\in\mathcal H(\delta,\Delta)\right\}$, we get that
$$
\dot v(t)\leq -k_{\delta,\Delta} v(t)\,,\quad \forall t\in\Gamma\,,
$$
which corresponds to \rref{V1dot6} with $c_3(s)\equiv 0$. The rest of the proof follows along the same lines as in Section \ref{prooftheo01} by noticing that both Assumption \ref{A23} and the bound on the gradient were only used in order to establish \rref{V1dot6}.


\subsection{Proof of Proposition \ref{propbound}}\label{sectionproofpropbound}
\vspace{-0.5cm}We claim that, whenever $V(t,x)=\underline\alpha(b)$, its derivative along the trajectories of \rref{sys1}, which we denote by $\dot V$, is non positive. To this end, notice that \rref{table1} implies that, if $V(t,x)=\overline\alpha(b)$, then $x\in\mathcal H(\overline\alpha^{-1}\circ\underline \alpha(b),b)$, which is nonempty (since $\underline\alpha(b)\leq\overline\alpha(b)$) and included in $\mathcal H(a,b)$ (since it is assumed that $\overline\alpha(a)<\underline\alpha(b)$). Hence, the claim is proved in view of \rref{table2}. For any $t_0\in\mRp$ and any $x_0\in \mR^n$, by defining $v(t):=V(t,x(t,t_0,x_0))$, we therefore get that, for all $t\geq t_0$,
$$
v(t)=\underline\alpha(b)\quad\Rightarrow\quad \dot v(t)\leq 0\,,
$$
which ensures in its turn, by the continuity of $v(t)$, that
$$
v(t_0)\leq \underline\alpha(b)\quad\Rightarrow\quad v(t)\leq\underline\alpha(b)\,,\qquad \forall t\geq t_0\,.
$$
The conclusion follows by noticing that, from \rref{table1},
$$
\norm{x_0}\leq \overline\alpha^{-1}\circ\underline\alpha(b)\quad \Rightarrow\quad v(t_0)\leq \underline\alpha(b)\,, \qquad\textrm{and}\qquad v(t)\leq\underline\alpha(b)\quad \Rightarrow\quad \norm{x(t)}\leq b\,.
$$


\section{Conclusion}\label{sec:concl}

We have presented results for uniform semiglobal practical asymptotic stability of nonlinear time-varying systems. Our main theorem establishes that the cascade of two USPAS systems remains USPAS; it relies on a condition of local boundedness of the solutions of the cascade and the knowledge of a Lyapunov function for the perturbed subsystem, in the absence of the interconnection term. As a corollary, we have that, under similar conditions, the cascade of two uniformly semiglobally asymptotically stable systems remains USAS. The latter generalizes, in its turn, previous results reported in the literature. Further research is carried out in the direction of extending these results to \emph{global} practical properties and on deriving sufficient conditions for boundedness of solutions. 

We have also illustrated the usefulness of our main results in the PID control of manipulators with external disturbances. The case-study presented here addresses, as far as we know, an important open problem in the literature of robot control. 


\small 
\bibliographystyle{plain}
\bibliography{refs}

\newcommand{\SortNoop}[1]{} 
  \def\nesic{Ne\v{s}i\'{c}\,} %
  \def\astrom{{\SortNoop{As}\AA}str{\"{o}}m\,}\let\c=\cedille
\begin{thebibliography}{10}

\bibitem{AILOGI97}
A.~Ailon, R.~Lozano-Leal, and M.~Gil'.
\newblock Point-to-point regulation of a robot with flexible joints including
  electrical effects of actuator dynamics.
\newblock {\em IEEE Trans. on Automat. Contr.}, 1997.

\bibitem{CHLOautomatica05}
A.~Chaillet and A.~Lor\'{\i}a.
\newblock A converse theorem for uniform semiglobal practical asymptotic
  stability: application to cascaded systems.
\newblock {\em Accepted in Automatica}, 2006.

\bibitem{THORB}
T.~I. Fossen.
\newblock {\em Guidance and control of ocean vehicles}.
\newblock {John Wiley \& Sons Ltd.}, 1994.
\newblock ISBN: 0-471-94113-1.

\bibitem{JIATEEPRA}
Z.~P. Jiang, A.~Teel, and L.~Praly.
\newblock Small gain theorems for {ISS} systems and applications.
\newblock {\em Math. of Cont. Sign. and Syst.}, 7:95--120, 1994.

\bibitem{KELLYSBOOK}
R.~Kelly, V.~Santibanez, and A.~\loria.
\newblock {\em Control of Robot Manipulators in Joint Space}.
\newblock Springer, {Berlin, Germany}, 2005.

\bibitem{KKK}
M.~Krsti\'c, I.~Kanellakopoulos, and P.~Kokotovi\'c.
\newblock {\em Nonlinear and Adaptive control design}.
\newblock {John Wiley \& Sons, Inc.}, New York, 1995.

\bibitem{ERJENSTHESIS}
A.~A.~J. Lefeber.
\newblock {\em Tracking control of nonlinear mechanical systems}.
\newblock PhD thesis, University of Twente, {Enschede, The Netherlands}, 2000.

\bibitem{LINSONWAN}
Y.~Lin, E.~D. Sontag, and Y.~Wang.
\newblock A smooth converse {Lyapunov} theorem for robust stability.
\newblock {\em SIAM J. on Contr. and Opt.}, 34:124--160, 1996.

\bibitem{CASCFAP}
A.~\loria\ and E.~Panteley.
\newblock {\em {\em Cascaded nonlinear time-varying systems: analysis and
  design}}, volume 244 of {\em Lecture Notes in Control and Information
  Sciences}, chapter in {\em New directions in nonlinear observer design}.
\newblock {Springer Verlag}, {F. Lamnabhi-Lagarrigue, A. \loria, E. Panteley,
  eds., London}, 2004.

\bibitem{DTCASCTAC}
D.~\nesic\ and A.~\loria.
\newblock On uniform asymptotic stability of time-varying parameterized
  discrete-time cascades.
\newblock {\em IEEE Trans. on Automat. Contr.}, 2004.

\bibitem{NESTEEPK}
D.~\nesic, A.~Teel, and P.~Kokotovi\'c.
\newblock Sufficient conditions for stabilization of sampled-data nonlinear
  systems via discrete-time approximations.
\newblock {\em Syst. \& Contr. Letters}, 38:259--270, 1999.

\bibitem{NESTEESON}
D.~\nesic, A.~Teel, and E.~Sontag.
\newblock Formulas relating {${\cal K}{\cal L}$} stability estimates of
  discrete-time and sampled-data nonlinear systems.
\newblock {\em Syst. \& Contr. Letters}, 38:49--60, 1999.

\bibitem{PI2D}
R.~Ortega, A.~\loria\, and R.~Kelly.
\newblock A semiglobally stable output feedback {PI$^2$D} regulator for robot
  manipulators.
\newblock {\em IEEE Trans. on Automat. Contr.}, 40(8):1432--1436, 1995.

\bibitem{PBCELS}
R.~Ortega, A.~\loria\, P.~J. Nicklasson, and H.~Sira-Ram\'{\i}rez.
\newblock {\em {Passivity-based Control of Euler-Lagrange Systems: Mechanical,
  Electrical and Electromechanical Applications}}.
\newblock Comunications and Control Engineering. Springer Verlag, London, 1998.
\newblock ISBN 1-85233-016-3.

\bibitem{CASCAUT}
E.~Panteley and A.~Lor\'{\i}a.
\newblock Growth rate conditions for stability of cascaded time-varying
  systems.
\newblock {\em Automatica}, 37(3):453--460, 2001.

\bibitem{PANORT}
E.~Panteley and R.~Ortega.
\newblock Cascaded control of feedback interconnected systems: {Application} to
  robots with {AC} drives.
\newblock {\em Automatica}, 33(11):1935--1947, 1997.

\bibitem{PRAWAN}
L.~Praly and Y.~Wang.
\newblock Stabilization in spite of matched unmodelled dynamics and an
  equivalent definition of input-to-state stability.
\newblock {\em Math. of Cont. Sign. and Syst.}, 9:1--33, 1996.

\bibitem{SCISIC}
L.~Sciavicco and B.~Siciliano.
\newblock {\em Modeling and control of robot manipulators}.
\newblock {McGraw Hill}, {New York}, 1996.

\bibitem{SEISUA}
P.~Seibert and R.~Su{\'a}rez.
\newblock Global stabilization of nonlinear cascaded systems.
\newblock {\em Syst. \& Contr. Letters}, 14:347--352, 1990.

\bibitem{SEPBOOK}
R.~Sepulchre, M.~Jankovi\'c, and P.~Kokotovi\'c.
\newblock {\em Constructive nonlinear control}.
\newblock {Springer Verlag}, 1997.

\bibitem{SON89}
E.~D. Sontag.
\newblock Remarks on stabilization and {Input-to-State} stability.
\newblock In {\em Proc. 28th. IEEE Conf. Decision Contr.}, pages 1376--1378,
  {Tampa, Fl}, 1989.

\bibitem{SONTAC03}
E.~D. Sontag.
\newblock A remark on the converging-input converging-state property.
\newblock {\em IEEE Trans. on Automat. Contr.}, 48, 2003.

\bibitem{SPOVID2}
M.~Spong, S.~Hutchinson, and M.~Vidyasagar.
\newblock {\em Robotics Modeling and Control}.
\newblock John Wiley \& Sons, New York, 2005.

\bibitem{SUSKOK}
H.~J. Sussman and P.~V. Kokotovi{\'c}.
\newblock The peaking phenomenon and the global stabilization of nonlinear
  systems.
\newblock {\em IEEE Trans. on Automat. Contr.}, 36(4):424--439, 1991.

\bibitem{TEPEAE98}
A.~Teel, J.~Peuteman, and D.~Aeyels.
\newblock Global asymptotic stability for the averaged implies semi-global
  practical stability for the actual.
\newblock In {\em Proc. 37th IEEE Conf. on Decision and Control}, Tempa,
  Florida, December 1998.

\bibitem{TEEPRA1}
A.~Teel and L.~Praly.
\newblock Global stabilizability and observability imply semiglobal
  stabilizability by output feedback.
\newblock {\em Syst. \& Contr. Letters}, 22:313--325, 1994.

\bibitem{TEEPRAconverse}
A.R. Teel and L.~Praly.
\newblock A smooth {Lyapunov} function from a class-{KL} estimate involving two
  positive semi-definite functions.
\newblock {\em ESAIM: COCV}, 5, 2000.

\end{thebibliography}

\end{document}